\documentclass[a4paper,12pt]{article}
\usepackage{amsfonts,amssymb, pdfpages, xcolor}

\title{The Cauchy problem for the standard one pressure system of two fluid flows with energy equations}
\author{Mathilde Colombeau,\\ \\
 Instituto de Matem\'atica, Estat\'istica e Ciencia da Computa\c c\~ao,\\Universidade  Estadual de Campinas, SP, Brazil.}

\begin{document}
\maketitle

\begin {abstract}
We construct  with full rigorous mathematical proof a family of approximate solutions to the Cauchy problem for the standard system of two fluid flows with energy equations and we pass to the limit by weak compactness to obtain Radon measures that satisfy the equations in a natural weak sense.  Our method provides  a  convergent numerical method for the numerical calculation of these Radon measures by reducing the system of partial differential equations in the case of these approximate solutions to a system of ordinary differential equations. 
We observe numerically on the standard Toumi shock tube problem that the Radon measures from our method agree  with the numerical solutions previously obtained by other authors  with various different numerical methods. In a subsequent  numerical paper, using a standard confident scheme with splittings and vanishing viscosity (independent on the above construction), we observe  exactly the numerical solution  given by our mathematical proof.  \end{abstract}

e-mail: m.colombeau@orange.fr\\
\textit{*this research has been partly done  thanks to  financial support of FAPESP, processo 2012/15780-9.}\\

\textbf{1.  Introduction}.\\ 
We consider the  standard system  used to model a mixture of two immiscible fluids from the conservation laws of mass, momentum and energy with the natural assumption that the pressures  are equal inside the two fluids, as derived in  \cite{Wendroff} p. 373 and solved numerically in \cite{Munk,FFMunk,FlattenRoe,EvjeFlatten}. We recall the system as stated in \cite{Wendroff} p. 373

\begin{equation}\frac{\partial}{\partial t}(\rho_i \alpha_i)+\frac{\partial}{\partial x}(\rho_i\alpha_i v_i)=0,\end{equation}
\begin{equation}\frac{\partial}{\partial t}(\rho_i \alpha_iv_i)+\frac{\partial}{\partial x}(\rho_i\alpha_i( v_i)^2)+\alpha_i\frac{\partial p}{\partial x}=0,\end{equation}
\begin{equation}\frac{\partial}{\partial t}(\rho_i\alpha_ie_i)+p\frac{\partial}{\partial t}\alpha_i +\frac{\partial}{\partial x}(\rho_i\alpha_ie_iv_i)+\frac{\partial}{\partial x}(p\alpha_iv_i)=0,   \end{equation}
 with the state laws \cite{Munk,FFMunk,FlattenRoe}
\begin{equation}p=(K_i-1)(\rho_ie_i-\frac{\rho_i (v_i)^2}{2})-K_ip^\infty_i,\ i=1,2,\end{equation}
where $\rho_i,v_i,e_i,\alpha_i$ are respectively the density, the velocity, the specific total energy and the volumic proportion of fluid $i$; $p$ is the unique pressure of the mixture, so that (4) provides a compatibility relation.  The constants $K_i$ and $p^\infty_i$ are  obtained from experimental measurements on the two fluids.\\


This system is  used in industry to model the gas kick phenomenon in offshore oil exploitation  \cite{Avelar,Bendiksen,Larsen} and to model the cooling in nuclear power stations   \cite{Bestion,WAHA3}. This system has various peculiarities which render its numerical study far more difficult than the classical case of one single fluid.  Since it has been observed \cite{Munk} p. 2620  that close numerical methods can give significantly different numerical solutions, discover which is the correct one is primarily important. \\

  The purpose of this paper is to  provide  with full rigorous mathematical proofs  that numerical solutions obtained in \cite{Munk,FFMunk,FlattenRoe} are  approximate solutions of the equations. We prove existence of a solution in a weak sense. Our method provides  a  convergent numerical method for the numerical calculation of these approximate solutions by reducing the system of partial differential equations in the case of these approximate solutions to a system of ordinary differential equations. To this end we first construct a family of smooth approximate solutions of the initial value problem and we prove that they tend to satisfy the equations in a weak sense, i.e. when plugged into the equations one proves that the result  tends to 0. Passage to the limit of approximate solutions by weak compactness gives radon measures in density, momentum and total energy that satisfy the equations in a natural weak sense.	Other authors \cite{Joseph1}-\cite{Joseph7},\cite{Kunzinger2} also introduced   sequences of approximate solutions to treat nonlinear PDEs. Of course, we do not study the very difficult problem of uniqueness, which seems  out of reach by the method in use here. \\

	Then, since our  approximate solutions are obtained as solutions of a system of six scalar ODEs in a Banach space, we calculate numerically these approximate solutions by means of faithful classical convergent numerical methods for ODEs.
	We observe numerically that one obtains    same results as  these authors  concerning step values and locations of the discontinuities   but we observe  very neat peaks in liquid velocity and in gas velocity, with well defined top values and widths, that appear in our  solution. Of course one could wonder the relevance of these peaks since they could be interpreted as  artefacts. To clarify this point we observe numerically the stability and evolution with time of these peaks. Further, in a purely numerical  paper (without proofs) \cite{Colombeaunumerique}, we obtain exactly the same result (i.e. the peaks) from a  standard confident numerical scheme with  splittings and vanishing viscosity. This last scheme can be adapted to both system (1-4) and to the system obtained in \cite{Munk} p.2595 by adding a supplementary term to ensure hyperbolicity. One observes that this additional term ensures the disappearance of these peaks.  The theoretical and numerical results  in this paper and the numerical results in \cite{Colombeaunumerique} show that these peaks are part of the solution of system (1-4) as it is, i.e. without additional terms. Further these peaks  are strikingly similar to  the experimental observations on the liquid flow rate $(m^3/s)$ when the gas  kick occurs \cite{Avelar} figure 4 there, the survey papers \cite{Avelarthese,Bezeira,Choi,Loiola} and references there. It had already been noticed by B. Keyfitz et al \cite{Keyfitz0}-\cite{KeyfitzSan} that some fully nonlinear systems which are non hyperbolic could nevertheless produce realistic physical solutions.\\


\textbf{2. Statement of  the main result.} The main result is the construction of approximate solutions. Then one passes easily to the limit by weak compactness and the convergence of the numerical method follows at once from the construction. We obtain the approximate solutions as solutions of a system of 6 scalar ODEs in the Banach space $\mathcal{C}(\mathbb{T})$ of all continuous functions on the torus $\mathbb{T}=\mathbb{R}/\mathbb{Z}$. First, a simplifying notation: if $V,v$ are real valued  functions of $(x,t,\epsilon)\in\mathbb{T}\times[0,T[\times]0,\eta[ , \ T>0,\eta>0$, we set\\

$[V,v]_x=-\frac{1}{\epsilon}\{V(x-\epsilon,t,\epsilon)v^+(x-\epsilon,t,\epsilon)-V(x,t,\epsilon)|v|(x,t,\epsilon)+$\begin{equation}V(x+\epsilon,t,\epsilon)v^-(x+\epsilon,t,\epsilon)\},\end{equation}
where $v^+=max(v,0)$ and $v^-=max(-v,0)$ so that $v=v^+-v^-$ and $|v|=v^++v^-$. The index $x$ in $[\dots]_x$  is used to avoid confusion with a simple parenthesis and because in 2-D we use  $[\dots]_y$. We introduce  two strictly positive functions $\epsilon\longmapsto\kappa_1(\epsilon)$  and  $\epsilon\longmapsto\kappa_2(\epsilon)$  that tend to 0 as fast as needed; they are used in the theoretical  proof. We also introduce a real number $ \lambda \in ]0,1[$ which is used in a convolution.\\

Physicists use small cubes of side $dx>0$ during times $t$ and $t+dt$ to state the laws of physics, then pass formally to the limit $dx=0$ and $dt=0$ to obtain the usual PDEs. Our method consists in letting $dt\rightarrow 0$ for fixed $dx=\epsilon>0$, then solve (theoretically from a mathematical proof  and numerically from convergent schemes for ODEs) for each value of $\epsilon>0$ the system of ODEs so obtained. We pass to the limit $\epsilon=0$  on the solutions of the ODEs.  For system (1-3) the method gives  existence of approximate solutions, then Radon measures at their limit by weak compactness and the method gives also their  numerical approximation. A physical justification of this method is given in section 5 below.\\ \\

 Now setting $r_i=\rho_i\alpha_i, \ i=1,2,$ we  state the ODEs we will use in the form
\begin{equation}\frac{d}{dt}r_i+[r_i,v_i]_x=\kappa_1(\epsilon),\end{equation}
\begin{equation}\frac{d}{dt}(r_iv_i)+[r_iv_i,v_i]_x+\alpha_i\partial_x \overline{p}=0,\end{equation}
\begin{equation}\frac{d}{dt}(r_ie_i)+[r_ie_i,v_i]_x+p\frac{d}{dt}\alpha_i+ (\partial_x \overline{p})\alpha_iv_i+p\partial_x(\overline{\alpha_iv_i})=\kappa_2(\epsilon),\end{equation}
 together with the algebraic equations
\begin{equation}p=(K_i-1)(\rho_ie_i-\frac{\rho_i (v_i)^2}{2})-K_ip_i^\infty,\end{equation}
\begin{equation}v_i=\frac{(r_iv_i)}{r_i}, \end{equation}
where $\overline{p} $ and $\overline{\alpha_iv_i}$ are respective mollifications of $p$ and $\alpha_iv_i$: 
\begin{equation}\overline{p}=p*(\phi)_{\epsilon^\lambda}, \ \overline{\alpha_iv_i}=(\alpha_iv_i)*(\phi)_{\epsilon^\lambda}, \ \ i=1,2,\end{equation}
 where $\phi\in\mathcal{C}_c^\infty(\mathbb{R}), \phi$ even, $\phi$ positive, $\int\phi(x)dx=1$ and $(\phi)_\mu(x)=\frac{1}{\mu}(\phi)(\frac{x}{\mu}).$\\
\\
Remark. Via the terms $[\dots]_x$, which provide some kind of semi-discretization, the ODEs (6-8) model a difference between space and time:  to model the space derivatives one cannot consider "`physical points"' having a length smaller than the edge of a cube that contains several thousand molecules, while time can be considered as made of mathematical points. The assumption that the function $\phi$ is even models the isotropy of space away from boundaries. In absence of  precisely defined  solutions for nonconservative systems we believe important to stick to physical intuition as much as possible (see also section 6 concerning molecular agitation) at the same time as to provide rigorous mathematical proofs.\\

We assume that the initial conditions satisfy  the requirements imposed by physics. More precisely $\forall x\in \mathbb{T} \ \  \rho_i^0(x)> 0$ and $ 0<\alpha_i^0(x)<1$ (to have a real mixture) with the compatibility condition $\alpha_1^0(x)+\alpha_2^0(x)=1, e_i^0(x)> 0, p^0(x)> 0,$ with the compatibility condition imposed by the state laws (9) and
$$ \rho_i^0\in L^1(\mathbb{T}), v_i^0\in L^\infty(\mathbb{T}), e_i^0\in L^\infty(\mathbb{T}).$$
 When $\epsilon\rightarrow 0$  we approximate these initial conditions, for instance by convolution as in (11), by smooth functions in $x$-variable, $\rho_i^0(x,\epsilon)>0$ and $ 0<\alpha_i^0(x,\epsilon)<1$, with the compatibility condition $\alpha_1^0(x,\epsilon)+\alpha_2^0(x,\epsilon)=1$, $e_i^0(x,\epsilon)>0, p^0(x,\epsilon)>0,$ with the compatibility condition (9) for all $\epsilon$  and such that there exists $const>0$ independent on $\epsilon$ such that
\begin{equation}\|\rho_i^0(.,\epsilon)\|_{L^1}\leq const, \|v_i^0(.,\epsilon)\|_{L^\infty}\leq const, \|e_i^0(.,\epsilon)\|_{L^\infty}\leq const.\end{equation} 
The strict inequalities $\rho_i^0(x,\epsilon)>0, \ 0<\alpha_i^0(x,\epsilon)<1, \ e_i^0(x,\epsilon)>0$ and $p_i^0(x,\epsilon)>0$ will be needed in the proof. \\ 

 We assume that the fluids are such that
\begin{equation} K_2>K_1>1\end{equation}
and
\begin{equation} K_2p_2^\infty>K_1p_1^\infty=0.\end{equation}
Indeed in the application to the Toumi shock tube problem in \cite{FFMunk,FlattenRoe,Munk} one has $K_1=1.4,K_2=2.8,p_1^\infty=0,p_2^\infty=8.5 .10^8.$  The fact that $K_1p_1^\infty=0$ plays a role in the proof of the theorem (to obtain  formulas (25-27)).\\

 We have to assume that the system of ODEs (6-11) does not lead to presence of a void region in any fluid. More precisely we have to assume that we consider solutions of (6-11) on $[0,T]$ which have the property that for $i=1,2$
\begin{equation}\exists  m>0 \hbox{ such that } r_i(x,t,\epsilon)=\alpha_i(x,t,\epsilon)\rho_i(x,t,\epsilon) \geq m>0 \ \forall x\in \mathbb{T},\forall t\in[0,T],\end{equation} $\forall \epsilon>0$  small enough, which means absence of void region in any fluid.\\

 
In the sequel we assume $T$ is arbitrarily large  to simplify the statement of the result. If not the result below holds on $[0,T[ $  only. This limitation of absence of a void region in any fluid is  needed in the proof and we have been unfortunately unable to prove the theorem from the equations and the initial data without assuming (15), which was possible  in \cite{ColombeauZeit, ColombeauJDE} in the particular case of state laws $p=f(\rho)$. One has checked numerically that (15) is  satisfied in the Toumi shock tube problem considered in this paper.\\

\textbf{Theorem .} \textit{Under the above assumption of absence of a void region  in any fluid and if $2\lambda<1$, the system of ODEs (6-11) has a global unique solution in positive time. This solution approximates system (1-4) in the following weak sense:  for all test function $\psi\in \mathcal{C}^\infty(\mathbb{T})$, the following limits hold $\forall t$ and for $i=1,2$ when} $\epsilon\rightarrow 0$
\begin{equation}\int_{\mathbb{T}} \{\frac{\partial}{\partial t}(\rho_i \alpha_i)(x,t,\epsilon)\psi(x)-(\rho_i\alpha_i v_i)(x,t,\epsilon)\} \partial_x\psi(x) \}dx\rightarrow 0,\end{equation}
\begin{equation} \int_{\mathbb{T}}\{\frac{\partial}{\partial t}(\rho_i \alpha_iv_i)(x,t,\epsilon)\psi(x)-(\rho_i\alpha_i v_i^2)(x,t,\epsilon) \partial_x\psi(x)+[\alpha_i(x,t,\epsilon)\frac{\partial}{\partial x}\overline{p}(x,t,\epsilon)]\psi(x)\}dx\rightarrow 0,\end{equation}

$$\int_{\mathbb{T}}\{\frac{\partial}{\partial t}(\rho_i\alpha_ie_i)(x,t,\epsilon)\psi(x) -(\rho_i\alpha_ie_iv_i)(x,t,\epsilon)\partial_x\psi(x)+p(x,t,\epsilon)\frac{\partial}{\partial t}(\alpha_i)(x,t,\epsilon)\psi(x)+$$\begin{equation}[\frac{\partial}{\partial x}\overline{p}(x,t,\epsilon)\alpha_iv_i(x,t,\epsilon)+p(x,t,\epsilon)\frac{\partial}{\partial x}(\overline{\alpha_iv_i})(x,t,\epsilon)]\psi(x)\}dx\rightarrow  0.   \end{equation}
\\
\textit{  The two state laws (4 or 9) are verified algebraically. Further  $\overline{p}$ and $\overline{\alpha_iv_i}$ are mollifications of $p$ and $\alpha_iv_i$  respectively and   when}  $\epsilon\rightarrow 0$
\begin{equation}\int_{\mathbb{T}}(\overline{p}-p)(x,t,\epsilon)\psi(x)dx\rightarrow 0\end{equation}  \textit{and}  
\begin{equation}\int_{\mathbb{T}}(\overline{\alpha_i v_i}-\alpha_iv_i)(x,t,\epsilon)\psi(x)dx\rightarrow 0,\end{equation}
\textit{for all test function} $\psi\in \mathcal{C}^\infty(\mathbb{T}).$
\\

 The terms containing these mollifications  originate from  the state law (9):  the state laws stem from experiments done at a macroscopic order of smallness which makes a great difference with the conservation laws. The  need of mollification  of physical variables involved in state laws appeared in \cite{ColombeauJMP,ColombeauZeit,ColombeauJDE}. \\

For fixed $\epsilon>0$ the physical variables $\overline{p}$ and $\overline{\alpha_iv_i}$ are of class $\mathcal{C}^\infty$ in $x,t$ variable, therefore all nonconservative products in the equations make sense for each fixed $\epsilon$. Our construction extends easily to more than two fluids (immediate) and in multidimension as in \cite{ColombeauZeit}.\\

Now, before the proof, we need some preparation.\\


\textbf{3. Preliminary calculations.}\\

$\bullet$\textit{Calculation of the volumic fractions for fixed $\epsilon$.}
The two state laws (9) valid for the same values of pressure imply
$$(K_1-1)(\rho_1e_1-\frac{\rho_1v_1^2}{2})-K_1p_1^\infty=(K_2-1)(\rho_2e_2-\frac{\rho_2v_2^2}{2})-K_2p_2^\infty.$$

Replacing  $\rho_i$ by $\frac{r_i}{\alpha_i}$, multiplying by $\alpha_1\alpha_2$ and setting $\alpha_1=\alpha,\alpha_2=1-\alpha$
 one obtains the following second order equation in unknown $\alpha$ in function of the 6 independent variables $r_i,r_iv_i$ and $r_ie_i$:

$$\alpha^2\{K_1p_1^\infty-K_2p_2^\infty\}+$$ $$\alpha\{-(K_1-1)(r_1e_1-\frac{r_1(v_1)^2}{2})-K_1p_1^\infty-(K_2-1)(r_2e_2-\frac{r_2(v_2)^2}{2})+K_2p_2^\infty\}+$$\begin{equation}\{(K_1-1)(r_1e_1-\frac{r_1(v_1)^2}{2})\}=0\end{equation}
\\
which, for each value $x,t,\epsilon$  has exactly one root $\alpha$ in the interval $]0,1[$ since we will prove that one has always $p(x,t,\epsilon)>0$: if $f(\alpha)$ denotes the left hand-side of (21) one has $f(0)=\alpha_1(p+K_1p_1^\infty)>0$ and $f(1)=-\alpha_2(p+K_2p_2^\infty)<0$.  This value $\alpha$ from (21) is a smooth function of  the independent variables $r_i, r_iv_i$ and $r_ie_i$.\\

$\bullet$\textit{Resolution of the energy equations in time derivative.} We develop $p\frac{d\alpha_j}{dt}$ in the form

\begin{equation}p\frac{d \alpha_j}{dt}=p\sum_{i=1}^2(\frac{d r_i}{dt}\frac{\partial \alpha_j}{\partial r_i}+\frac{d (r_iv_i)}{dt}\frac{\partial \alpha_j}{\partial (r_iv_i)}+\frac{d (r_ie_i)}{dt}\frac{\partial \alpha_j}{\partial (r_ie_i)}).\end{equation}

Using (22) the system of ODEs (6-8)  can be written in the form $M\frac{dX}{dt}=N$ where $M$ is a $6\times 6$ matrix and $X=(r_1,r_2,r_1v_1,r_2v_2,r_1e_1,r_2e_2)^t$. The two energy equations (8) appear respectively in the form
$$(1+p\frac{\partial\alpha_1}{\partial (r_1e_1)})\frac{d}{dt}(r_1e_1)+p\frac{\partial\alpha_1}{\partial (r_2e_2)}\frac{d}{dt}(r_2e_2)=f_1,$$ 
$$p\frac{\partial\alpha_2}{\partial (r_1e_1)}\frac{d}{dt}(r_1e_1)+(1+p\frac{\partial\alpha_2}{\partial (r_2e_2)})\frac{d}{dt}(r_2e_2)=f_2,$$ 
where $f_1$ and $f_2$ are functions involving terms without time derivatives once one uses (6,7) to replace $\frac{dr_i}{dt}$ and $\frac{d(r_iv_i)}{dt}$ from (22) by terms without time derivatives.
Therefore the  determinant of $M$ is  $1-p(\frac{\partial \alpha}{\partial(r_2e_2)}-\frac{\partial \alpha}{\partial(r_1e_1)})$.
The system of ODEs can be resolved in time when this determinant is nonzero. So we compute this determinant. Differentiation in the variable $(r_1e_1)$ of  (21) with $\alpha=\alpha(r_1,r_2,r_1v_1,r_2v_2,r_1e_1,r_2e_2)$ gives
$$\frac{\partial \alpha}{\partial(r_1e_1)}[2\alpha(K_1p_1^\infty-K_2p_2^\infty)+\{\{\dots\}\}]=\alpha(K_1-1)-(K_1-1)$$
where we denote by $\{\{\dots\}\}$ the coefficient of $\alpha$ in (21). Similarly
$$\frac{\partial \alpha}{\partial(r_2e_2)}[2\alpha(K_1p_1^\infty-K_2p_2^\infty)+\{\{\dots\}\}]=\alpha(K_2-1).$$
One obtains that
\begin{equation} p(\frac{\partial \alpha}{\partial(r_2e_2)}-\frac{\partial \alpha}{\partial(r_1e_1)})=\frac{1}{2}\frac{[\alpha(K_2-K_1)+(K_1-1)]p}{\alpha(K_1p_1^\infty-K_2p_2^\infty)-p-K_1p_1^\infty}.\end{equation}
We have to check that this value is always different from 1 so that the system of the two energy ODEs (8) could be solvable in time derivative. Since $K_2-K_1>0, K_1-1>0$ (13), $0<\alpha<1$ and  since we will prove that $p(x,t,\epsilon)>0 \ \forall x,t,\epsilon$  the second member of (23) is different from 1. Therefore the determinant of $M$ is always nonzero .\\

Therefore, after multiplication by $M^{-1}$, for each fixed $\epsilon$ one has a system of 6 ODEs in the Banach space $\mathcal{C}(\mathbb{T})$ of all continuous functions on $\mathbb{T}$ of the form 
\begin{equation}\frac{dX}{dt}=F(X)\end{equation}
with $F:\Omega\subset(\mathcal{C}(\mathbb{T}))^6\longmapsto (\mathcal{C}(\mathbb{T}))^6$ where $\Omega=\{(r_1,r_2) \ / \ r_i>0,i=1,2\}\times(\mathcal{C}(\mathbb{T}))^4$, because of division by $r_i$ in (10). For fixed $\epsilon>0$ $F$  has Lipschitz coefficients uniform in $\Omega_M:=\Omega\cap\{inf\{r_i\}\geq \frac{1}{M}, \|r_i\|_\infty\leq M, \|r_iv_i\|_\infty\leq M,\|r_ie_i\|_\infty\leq M, \ i=1,2,\}$ for any $M>0$ ( the $\partial_x$ derivatives in (7,8) are done after convolution (11)). This remark will play a basic role in the sequel under the form that as long as a solution $(r_i,r_iv_i,r_ie_i), \ i=1,2,$ defined on some interval $[0,\delta(\epsilon)[$ takes its values in some set $\Omega_M$ then this solution can be extended to a larger interval $[0,\delta(\epsilon)+\eta(\epsilon)[, \ \eta(\epsilon)>0$. This will permit to prove from suitable a priori estimates existence of a global solution of (6-11) in positive time for each fixed $\epsilon$.\\

$\bullet$\textit{An ODE satisfied by the pressure.} From the state laws (9) for the fluid 1, in  which one uses that $K_1p_1^\infty=0$, one has
$$\partial_t(\rho_1\alpha_1e_1)=\frac{\partial_t(p\alpha_1)}{K_1-1}+\frac{\partial_t(\rho_1\alpha_1(v_1)^2)}{2}$$
and
$$[\rho_1\alpha_1e_1,v_1]_x=\frac{[p\alpha_1,v_1]_x}{K_1-1}+\frac{[\rho_1\alpha_1(v_1)^2,v_1]_x}{2}.$$
Plugging these two formulas into the energy equation (8) with $i=1$ and recalling that $r_1=\rho_1\alpha_1$ one obtains

$$\frac{\partial_t(p\alpha_1)}{K_1-1}+\frac{\partial_t(r_1(v_1)^2)}{2}+\frac{[p\alpha_1,v_1]_x}{K_1-1}+\frac{[r_1(v_1)^2,v_1]_x}{2}+p\partial_t\alpha_1+(\partial_x\overline{p})\alpha_1v_1+p\partial_x(\overline{\alpha_1v_1})=\kappa_2(\epsilon)$$
i.e.\\
\\
$\partial_t(p\alpha_1)+[p\alpha_1,v_1]_x+\frac{K_1-1}{2}\{\partial_t(r_1(v_1)^2)+[r_1(v_1)^2,v_1]_x+2(\partial_x\overline{p})\alpha_1v_1+2p\partial_x(\overline{\alpha_1v_1})\}=$\begin{equation}-(K_1-1)p\partial_t\alpha_1+(K_1-1)\kappa_2(\epsilon).\end{equation}

Auxiliary calculation: to the Euler equation (7) for $i=1$ we subtract the continuity equation (6) for $i=1$ multiplied by $v_1$. This calculation is justified because the solutions for fixed $\epsilon>0$ are smooth (one can regularize $v_i^\pm$ by setting $v_i^+-v_i^-=v_i,v_i^++v_i^-=\sqrt{v_i^2+a^2}, a\in\mathbb{R}, a>0$; nothing is changed in the results and proofs see \cite{ColombeauZeit} p. 2586-2587). One obtains
$$r_1\partial_tv_1=[r_1,v_1]_xv_1-v_1\kappa_1(\epsilon)-[(r_1v_1),v_1]_x-\alpha_1\partial_x\overline{p}.$$
Using $\partial_t(r_1(v_1)^2)=r_1v_1\partial_tv_1+v_1\partial_t(r_1v_1)$, the above formula and (7), the term $\frac{K_1-1}{2}\{\dots\}$ in (25) becomes\\
\\
$\frac{K_1-1}{2}\{[r_1,v_1]_x(v_1)^2-(v_1)^2\kappa_1(\epsilon)-[r_1v_1,v_1]_xv_1-\alpha_1(\partial_x\overline{p})v_1-$ \begin{equation}[r_1v_1,v_1]_xv_1-\alpha_1(\partial_x\overline{p})v_1+[r_1(v_1)^2,v_1]_x+2(\partial_x\overline{p})\alpha_1v_1+2p\partial_x(\overline{\alpha_1v_1})\}.\end{equation}
\\
Developping from (5) the  brackets $[\dots]_x$ that are in (26) one obtains \\ 
\\
$$\{[r_1,v_1]_x(v_1)^2-2[r_1v_1,v_1]_xv_1+[r_1(v_1)^2,v_1]_x\}(x,t,\epsilon)=\frac{1}{\epsilon}\{-r_1v_1^+(x-\epsilon,t,\epsilon)(v_1(x-\epsilon,t,\epsilon)-$$ $$v_1(x,t,\epsilon))^2-r_1v_1^-(x+\epsilon,t,\epsilon)(v_1(x+\epsilon,t,\epsilon)-v_1(x,t,\epsilon))^2\}.$$
\\
Plugging this last result into (26), then (25)  one finally obtains by developing $-[p\alpha_1,v_1]_x$ according to (5)\\
\\
$$\partial_t(\alpha_1p)(x,t,\epsilon)=\{\frac{1}{\epsilon}\{\alpha_1pv_1^+(x-\epsilon,t,\epsilon)-\alpha_1p|v_1|(x,t,\epsilon)+\alpha_1pv_1^-(x+\epsilon,t,\epsilon)\}+$$ $$\frac{K_1-1}{2\epsilon}\{r_1v_1^+(x-\epsilon,t,\epsilon)(v_1(x-\epsilon,t,\epsilon)-v_1(x,t,\epsilon))^2+r_1v_1^-(x+\epsilon,t,\epsilon)(v_1(x+\epsilon,$$
$t,\epsilon)-v_1(x,t,\epsilon))^2\}-(K_1-1)p\partial_x(\overline{\alpha_1v_1})(x,t,\epsilon)+\frac{K_1-1}{2}(v_1)^2(x,t,\epsilon)\kappa_1(\epsilon)-$ \begin{equation}(K_1-1)p\partial_t\alpha_1(x,t,\epsilon)+(K_1-1)\kappa_2(\epsilon)\}.\end{equation}

This formula will be used in the proof by noting that the terms in second member are positive except possibly the 3 terms involving the factor $p(x,t,\epsilon)$ since $r_1\geq 0,\alpha_1\geq 0,p\geq 0,K_1> 1,v_1^\pm\geq 0$. In the sequel of the  proof we will consider values of $(x,t)$ such that $p(x,t,\epsilon)=0$,  therefore, for these values,  one will have $\partial_t (\alpha_1 p)(x,t,\epsilon)>0$ since $\kappa_2(\epsilon)>0$ and all other terms are positive or null. \\


$\bullet$ \textit{A simplification in an integral.}
\begin{equation}\int_{\mathbb{T}}\{\alpha_iv_i(\partial_x\overline{p})+p\partial_x(\overline{\alpha_iv_i})\}(x,t,\epsilon)dx=0.\end{equation}

proof. Dropping the notations $t$ and $\epsilon$, from (11) this integral is equal to\\

$\int\alpha_iv_i(x)p(x-\epsilon^\lambda\mu)\frac{1}{\epsilon^\lambda}\phi'(\mu)d\mu dx+\int\alpha_iv_i(x-\epsilon^\lambda\mu)p(x)\frac{1}{\epsilon^\lambda}(\phi)'(\mu)d\mu dx=$
$$\int\alpha_iv_i(x)p(x-\epsilon^\lambda\mu)\frac{1}{\epsilon^\lambda}(\phi'(\mu)+\phi'(-\mu))d\mu dx.$$ Then use the assumption that $\phi$ is even.$\Box$\\

Now one can start the proof of the theorem by a priori estimates.\\


\textbf{4. Proof of the theorem.} From the local existence-uniqueness theorem for ODEs in the Lipschitz case, for all $\epsilon>0$ there exists $\delta(\epsilon)>0$, depending on $\epsilon$, and a unique solution of  system (6-11) on $[0,\delta(\epsilon)[$. From the choice of approximations of the initial conditions (strict positiveness) we choose $\delta(\epsilon)$ small enough so that $\forall t\in[0,\delta(\epsilon)[, \ \forall x\in\mathbb{T}, \forall i=1,2$
\begin{equation}\rho_i(x,t,\epsilon)>0, \ p(x,t,\epsilon)>0, \ e_i(x,t,\epsilon)>0, \ 0<\alpha_i(x,t,\epsilon)<1.\end{equation}

Now one assumes that the solution exists on some non necessarily small interval $[0,\delta(\epsilon)[$, for a given finite value $\delta(\epsilon)$, and that on $[0,\delta(\epsilon)[$ this solution satisfies (29). Our aim is to obtain, under assumption (15) of absence of void region in one fluid, $L^1$ and $L^\infty$ a priori estimates on the solution on $[0,\delta(\epsilon)[$, some of them uniform in $\epsilon$, so as to extend this solution to the right of $\delta(\epsilon)$ and also to use them in the proof of the approximations (16-20). The notation $const$ will denote  values independent on $\epsilon$ and on $t\in [0,\delta(\epsilon)[$.\\
\\
$\bullet$ \textit{First step: $L^1$ bounds uniform in $\epsilon$.}
From (6,5) and replacement of $|v|$ by $v^++v^-$ we obtain 
$$\frac{d}{dt}r_i(x,t,\epsilon)=\frac{1}{\epsilon}[r_iv_i^+(x-\epsilon,t,\epsilon)-r_iv_i^+(x,t,\epsilon)-r_iv_i^-(x,t,\epsilon)+r_iv_i^-(x+\epsilon,t,\epsilon)]+\kappa_1(\epsilon).$$
The two integrals in $v_i^+$ simplify by translation when one integrates in $x$ on $\mathbb{T}$, as well as the two integrals in $v_i^-$. One obtains
$$\frac{d}{dt}\int r_i(x,t,\epsilon)dx=\kappa_1(\epsilon) \ \forall t\in [0,\delta(\epsilon)[,$$
i.e. since $r_i=\rho_i\alpha_i>0$ (29), $ \kappa_1(\epsilon)\rightarrow 0$  and $\delta(\epsilon)$ finite
\begin{equation}\|r_i(.,t,\epsilon)\|_{L^1}\leq const,\end{equation}
where $L^1$ denotes $L^1(\mathbb{T})$. We add the two energy equations (8) so as to eliminate the nonconservative terms $p\frac{d}{dt}\alpha_i$ since $\alpha_1+\alpha_2=1$. Then in the same way as for (30), using (5) and (28), one obtains a  $L^1$ bound for $r_1e_1+r_2e_2$ uniform in $\epsilon$, from the positiveness of $r_ie_i=\rho_i\alpha_ie_i$ (29). Therefore 
\begin{equation}\|r_ie_i(.,t,\epsilon)\|_{L^1}\leq const, \ i=1,2  \ \forall t\in [0,\delta(\epsilon)[,\end{equation}
which implies from the state laws (9) and $r_i=\rho_i \alpha_i, 0<\alpha_i<1$, that
\begin{equation} \|r_iv_i^2(.,t,\epsilon)\|_{L^1}\leq const, \ i=1,2 \ \forall t\in [0,\delta(\epsilon)[ ,\end{equation} since $K_ip_i^\infty\geq 0$ and $\alpha_ip>0$ from (29).
The state laws (9)  and the  positiveness of  $\alpha_i p$ imply also from (31) that
\begin{equation}\|(\alpha_i p)(.,t,\epsilon)\|_{L^1}\leq const, \ i=1,2 \ \forall t\in [0,\delta(\epsilon)[.\end{equation}
Since $\alpha_1+\alpha_2=1$ (33) implies 
\begin{equation}\| p(.,t,\epsilon)\|_{L^1}\leq const \ \forall t\in [0,\delta(\epsilon)[.\end{equation}
Since $|v_i|\leq max(1,v_i^2)$ (30,32) imply
\begin{equation}\|(r_iv_i)(.,t,\epsilon)\|_{L^1}\leq const, \ i=1,2 \ \forall t\in [0,\delta(\epsilon)[.\end{equation}


After these $L^1$ bounds uniform in $\epsilon$ that followed  easily from the equations (6-11) and from the positiveness assumption (29), we are going to obtain  $L^\infty$ bounds depending on $\epsilon$ that will permit to prove the existence of a global solution to the ODEs (6-11) for fixed $\epsilon$.\\

$\bullet$ \textit{Second step: $L^\infty$ bounds depending on $\epsilon$.}
The bound (34) and the convolutions (11) imply 
\begin{equation}\|\overline{p}(.,t,\epsilon)\|_{L^\infty}\leq \frac{const}{\epsilon^\lambda}, \ \|(\partial_x\overline{p})(.,t,\epsilon)\|_{L^\infty}\leq \frac{const}{\epsilon^{2\lambda}} \ \forall t\in [0,\delta(\epsilon)[.\end{equation}
The basic point lies in bounds of $\|v_i\|_\infty$, which is more delicate. From (36) a proof based on (6,7) similar to the one in \cite{ColombeauJDE} pp. 204-205 to prove (31) there or to the more detailed one in 
\cite{ColombeauZeit} pp. 2581-2582, formulas (27-29,24) there, gives
\begin{equation}\|v_i(.,t,\epsilon)\|_\infty\leq \frac{const}{\epsilon^{2\lambda}},\ i=1,2\ \forall t\in [0,\delta(\epsilon)[,\end{equation}
assuming absence of void regions (15) and small enough values $\epsilon>0$.
 This proof is sketched as follows.  From (5,7) one develops at order one in $dt$ the quantity $(r_iv_i)(x,t+dt,\epsilon)$ at time $t$ for an arbitrarily small positive increment $dt$, as well as $r_i(x,t+dt,\epsilon)$ from (5,6). Then one uses these developments to develop  the quotient $v_i(x,t+dt,\epsilon)=\frac{(r_iv_i)(x,t+dt,\epsilon)}{r_i(x,t+dt,\epsilon)}$. The first term we obtain for  the quotient is bounded by a barycentric combination bounded by $\|v_i(.,t,\epsilon)\|_\infty$. The second term  is bounded by $\frac{\|\partial_x\overline{p}\|_\infty}{m}$ from (15). One obtains 
$$|v_i(x,t+dt,\epsilon)|\leq\|v_i(.,t,\epsilon)\|_\infty+dt\frac{\|\partial_x\overline{p}\|_\infty}{m}+dt \ R(x,t,\epsilon,dt)$$
where $R$ is a remainder obtained from the remainders in 
 the first order developments of  $(r_iv_i)(x,t+dt,\epsilon)$  and $r_i(x,t+dt,\epsilon)$. This remainder $R$  disappears finally by dividing the interval $[0,t]$ into small intervals of length $dt=\frac{t}{n}$ and letting $n\rightarrow +\infty$, thus giving (37). In this proof we used assumption (15) of absence of void region in the fluid $i$. \\

We prefer to give a proof of (38) below simply based on (37) not to use assumption (15).\\
\\
From (6,5) since $r_i$ and $v_i^\pm$ are positive
$$\frac{d r_i}{dt}(x,t,\epsilon)\geq -\frac{1}{\epsilon}r_i(x,t,\epsilon)|v_i(x,t,\epsilon)|+\kappa_1(\epsilon).  $$
 From (37) and $r_i(x,0,\epsilon)\geq m>0$ an explicit solution of the elementary ODE $X'(t)=-AX(t)+B(t)$ implies that 
\begin{equation}\exists m_1(\epsilon)>0 \hbox{ such that } r_i(x,t,\epsilon)\geq m_1(\epsilon) \ \forall t\in [0,\delta(\epsilon)[, \end{equation}
  for some value $m_1(\epsilon)>0$, see \cite{ColombeauZeit} pp. 2584,2585. Since $|\alpha_iv_i|\leq |v_i|$, (37) and the convolution (11) imply
\begin{equation}\|\overline{\alpha_iv_i}(.,t,\epsilon)\|_\infty\leq M(\epsilon), \ \|\partial_x(\overline{\alpha_iv_i})(.,t,\epsilon)\|_\infty\leq M(\epsilon)\end{equation}  $\forall t\in [0,\delta(\epsilon)[$, for some values $M(\epsilon)$ independent on $t\in[0,\delta(\epsilon)[$.\\
\\
Now we obtain a bound of $\|r_i(.,t,\epsilon)\|_\infty$ from (6,5) and the $L^\infty$ bound (37) of $v_i$ as follows: 
\begin{equation}r_i(x,t,\epsilon)\leq r_i^0(x,\epsilon)+\int_0^t\frac{2}{\epsilon}\|r_i(.,\tau,\epsilon)\|_\infty\frac{const}{\epsilon^{2\lambda}}d\tau+\kappa_1(\epsilon) t.\end{equation} Gronwall formula implies
\begin{equation}\|r_i(.,t,\epsilon)\|_\infty\leq (\|r_i^0(.,\epsilon)\|_\infty+\kappa_1(\epsilon)\delta(\epsilon))exp(\frac{const}{\epsilon^{2\lambda+1}})\ \forall t\in [0,\delta(\epsilon)[.\end{equation}

The bounds (37) and (41) give  a bound 
\begin{equation}\|r_iv_i(.,t,\epsilon)\|_\infty\leq M(\epsilon) \ \forall t\in [0,\delta(\epsilon)[, \end{equation}
for some $M(\epsilon)$. 
Using the sum of the two energy equations (8), the $L^\infty$ bound (37) of $v_i$ for the second term in (8), the $L^\infty$ bounds (37,39) for $\alpha_iv_i$ and $\partial_x(\overline{\alpha_i v_i})$, the $L^\infty$ bounds (36) for  $\partial_x\overline{p}$ and the bound  $p\leq const \ r_ie_i$, i.e. a linear bound in $r_ie_i$,  one obtains for $r_1e_1+r_2e_2$ a bound similar to the bound (40) obtained for $r_i$. Then, from Gronwall formula, one obtains the following $L^\infty$ bound as (41):


$$\|(r_1e_1+r_2e_2)(.,t,\epsilon)\|_\infty\leq M(\epsilon)\ \forall t\in [0,\delta(\epsilon)[, $$
for some $M(\epsilon)$, which implies from (29) that
\begin{equation}\|(r_ie_i)(.,t,\epsilon)\|_\infty\leq M(\epsilon)\ \forall t\in [0,\delta(\epsilon)[.\end{equation}

The minoration (38) and the $L^\infty$ bounds (41,42,43) show  that on $[0,\delta(\epsilon)[$ the solution of the system of ODEs (6,11) for fixed $\epsilon$ remains in some set $\Omega_{M}\subset(\mathcal{C}(\mathbb{T}))^6$. Therefore as noticed above the solution extends to $[0,\delta(\epsilon)+\eta(\epsilon)[$ for some $\eta(\epsilon)>0$. Now we will prove that the strict positiveness property (29) is true on $[\delta(\epsilon),\delta(\epsilon)+\eta(\epsilon)[$, which will permit  a proof by absurd that the solution exists on $[0,+\infty[$.\\


$\bullet$ \textit{Third step: strict positiveness (29) on $[\delta(\epsilon),\delta(\epsilon)+\eta(\epsilon)[$.}
We do not  use directly assumption (15)  because the following proof will also be needed in pressure. Let us prove that, for fixed $i$ and fixed $\epsilon$,  $r_i(x,t,\epsilon)>0 \ \forall x\in \mathbb{T} \forall t\in [\delta(\epsilon),\delta(\epsilon)+\eta(\epsilon)[$.
By absurd let us assume that $\exists t_0\leq \delta(\epsilon)$ and $\exists x_0\in\mathbb{T}$ such that $r_i(x_0,t_0,\epsilon)=0$. We can choose them such that $r_i(x,t,\epsilon)>0 \ \forall x\in \mathbb{T} \forall t<t_0$. Obviously $\partial_t r_i(x_0,t_0,\epsilon)\leq 0$. Now (6) gives
 \begin{equation}\partial_t r_i(x_0,t_0,\epsilon)=\frac{1}{\epsilon}\{r_iv_i^+(x_0-\epsilon,t_0,\epsilon)-0+r_iv_i^-(x_0+\epsilon,t_0,\epsilon)\} +\kappa_1(\epsilon)\end{equation}
which shows that $\partial_t r_i(x_0,t_0,\epsilon)\geq \kappa_1(\epsilon)>0$. Therefore we obtain a contradiction. We have proved by absurd that $r_i(x,t,\epsilon)>0  \  \forall x\in \mathbb{T}  \  \forall t\in [\delta(\epsilon),\delta(\epsilon)+\eta(\epsilon)[$.  From the second degree equation (21) one has  $\alpha_i(x,t,\epsilon)>0$ since the equation has exactly one root in $]0,1[$.  From formula (27) the same proof by absurd applies to $\alpha_1 p$ since $\alpha_1>0$: indeed since $p(x_0,t_0,\epsilon)=0$ all 3 possibly negative terms in the right hand-side of (27) disappear and the same reasoning applies. Therefore  $p(x,t,\epsilon)>0 \  \forall x\in \mathbb{T}  \  \forall t\in [\delta(\epsilon),\delta(\epsilon)+\eta(\epsilon)[$; then the state laws (9) give the strict positiveness of $r_ie_i$. \\

$\bullet$ \textit{Fourth step: global solution of the ODEs.}
Now one can prove that for all $\epsilon>0$ the system of ODEs (6-11) has a global solution in positive time. By absurd let us assume the solution  with strict positiveness property (29) ceases to exist at some time $\delta(\epsilon)>0$. We proved above from (29) on $[0,\delta(\epsilon)[$ that there exists $\eta(\epsilon)>0$ such that the solution  can be continued on $[\delta(\epsilon),\delta(\epsilon)+\eta(\epsilon)[$ with again the strict positiveness property (29), which produces a contradiction.\\

Now it remains to prove that the solution of the system of ODEs satisfies the approximations (16-20).\\


$\bullet$ \textit{Fifth step: proof of the approximations (16-20).} The proofs of (16),(17) and (18) are similar. Let us consider (17) and let $I$ denote the left hand-side of (17). One replaces $\partial_t (\rho_i\alpha_iv_i)=\partial_t(r_iv_i)$ by its value using (7,5). This gives $$I=\int_{\mathbb{T}}\{\frac{1}{\epsilon}\{r_iv_iv_i^+(x-\epsilon,t,\epsilon)-r_iv_iv_i^+(x,t,\epsilon)-r_iv_iv_i^-(x,t,\epsilon)+r_iv_iv_i^-(x+\epsilon,t,\epsilon)\}\psi(x)$$ $-(\alpha_i\partial_x\overline{p})(x,t,\epsilon)\psi(x)-(r_i(v_i)^2)(x,t,\epsilon)\psi'(x)+(\alpha_i\partial_x\overline{p})(x,t,\epsilon)\psi(x)\}dx.$
\\
\\
The two terms $\alpha_i\partial_x\overline{p}$ simplify. One replaces $r_iv_i^2$ by $r_iv_i(v_i^+-v_i^-)$. Then $I$ is the sum of a term  involving $v_i^+$ and a similar term involving $v_i^-$. The term involving $v_i^+$ is 
$$\int_{\mathbb{T}}\{\frac{1}{\epsilon}[r_iv_iv_i^+(x-\epsilon,t,\epsilon)-r_iv_iv_i^+(x,t,\epsilon)]\psi(x)-r_iv_iv_i^+(x,t,\epsilon)\psi'(x)\}dx.$$
After a change of variable this integral becomes
\begin{equation}\int_{\mathbb{T}}\{r_iv_iv_i^+(x,t,\epsilon)[\frac{\psi(x+\epsilon)-\psi(x)}{\epsilon}-\psi'(x)]\}dx= \int_{\mathbb{T}} r_iv_iv_i^+(x,t,\epsilon)O(\epsilon)dx.\end{equation}
Using $\|r_iv_i(.,t,\epsilon)\|_{L^1}\leq const$ (35) and $\|v_i^+\|_\infty\leq \frac{const}{\epsilon^{2\lambda}}$
(37) one obtains that the  integral in (45) is equal to $O(\epsilon^{1-2\lambda})$. Same result for the term involving $v_i^-$ and finally $I=O(\epsilon^{1-2\lambda})$. This proof applies as well to (16) and (18) (for (18) the two nonconservative terms $p\frac{\partial}{\partial t}\alpha_i$ simplify). The proof of (19) is immediate from (34). For (20) one uses Holder's formula,  assumption (15) on  absence of void region and (32): \begin{equation}\int |\alpha_iv_i|dx\leq\int |v_i|dx\leq (\int(v_i)^2dx)^{\frac{1}{2}}\leq(\int\frac{r_i}{m} (v_i)^2dx)^{\frac{1}{2}}\leq const.\end{equation}$\Box$\\

\textbf{5. Passage to the limit in approximate solutions.} We proved that the families $(\alpha_{i,\epsilon}\rho_{i,\epsilon}), (\alpha_{i,\epsilon}\rho_{i,\epsilon}e_{i,\epsilon}), (\alpha_{i,\epsilon}\rho_{i,\epsilon}v_{i,\epsilon}), (\alpha_{i,\epsilon}\rho_{i,\epsilon}v_{i,\epsilon}^2)$ and $(p_{i,\epsilon})$ are $L^1$-stable (30-35). Therefore they are relatively compact for the *weak topology of the space $M(\mathbb{T}\times ]0,T[)$ of Radon measures on $\mathbb{T}\times ]0,T[$ \cite{Brezis}. Consequently there are convergent subsequences to Radon measures denoted respectively $(\alpha_i\rho_i), (\alpha_i\rho_ie_i), (\alpha_i\rho_iv_i), (\alpha_i\rho_iv_i^2)$ and $p$.\\

 If $\psi\in\mathcal{C}_c^\infty(\mathbb{R}\times]0,T[)$ it follows from (16) that
$$\partial_t(\rho_i\alpha_i)+\partial_x(\rho_i\alpha_iv_i)=0$$ in the sense of distributions. The momentum equations are slightly more complicated because of the nonconservative terms $\alpha_i\partial_xp$. Since $\rho_{i,\epsilon}\alpha_{i,\epsilon}v_{i,\epsilon}$ and $\rho_{i,\epsilon}\alpha_{i,\epsilon}v_{i,\epsilon}^2 $ tend in the  *weak topology (for a subsequence) to Radon measures it follows from (17) that
 $$\int \alpha_{i,\epsilon}(x,t)\partial_x\overline{p}_\epsilon(x,t)\psi(x,t)dxdt\rightarrow <-\partial_t(\rho_{i}\alpha_{i}v_{i})+\partial_x(\rho_{i}\alpha_{i}v_{i}^2),\psi>$$
where derivations and brackets are intended in the sense of distributions. Therefore $\alpha_{i,\epsilon}\partial_x\overline{p}_\epsilon$ tends to a distribution denoted $\alpha_i\partial_x\overline{p}=-\partial_t(\rho_{i}\alpha_{i}v_{i})+\partial_x(\rho_{i}\alpha_{i}v_{i}^2)$ and the momentum equations are satisfied in the sense of distributions.\\

 Now if further we assume that the velocities $v_{i,\epsilon}$ are $L^\infty$-stable (as physically needed but that we could not prove) the families $(\rho_{i,\epsilon}\alpha_{i,\epsilon}e_{i,\epsilon}v_{i,\epsilon}), i=1,2$ are $L^1$-stable. One can pass to the limit in a subsequence in the space of Radon measures for the *weak topology which gives a Radon measure denoted $\rho_i\alpha_ie_iv_i$. From (21) $\alpha_{i,\epsilon}$ is differentiable in $t$ valued in $L^\infty(\mathbb{T})$ and $\partial_t\alpha_i\in L^\infty(\mathbb{T})$ therefore the family $p_\epsilon\partial_t\alpha_{i,\epsilon}$is $L^1$-stable uniformly in $t$. Let $p\partial_t\alpha_i\in M(\mathbb{T})\times]0,T[$ be its limit (for the same subsequence as above for all the other limits). One proves that the term $(\partial_x\overline{p}_\epsilon\alpha_{i,\epsilon}v_{i,\epsilon}+p_\epsilon\partial_x(\overline{\alpha_{i,\epsilon}v_{i,\epsilon}})$  tends to a distribution from (18) since the first three terms in (18) tend to distributions. Then the energy equations are satisfied in the sense of distributions in the form $$\partial_t(\rho_i\alpha_ie_i)+\partial_x(\rho_i\alpha_ie_iv_i)+(p\partial_t\alpha_i)+(\partial_x\overline{p}\alpha_{i}v_{i}+p_\epsilon\partial_x(\overline{\alpha_{i,\epsilon}v_{i,\epsilon}}))=0$$ in which we recall that each of the four terms is defined as a distribution by weak compactness.\\

If $\psi\in\mathcal{C}_c^\infty(\mathbb{T}\times[0,T[)$ is such that $\psi(.,0)\not=0$ one can state the exact equations taking into account the initial conditions as usual from integration by parts in time for $\epsilon>0$ then passage to the limit.\\

These Radon measures  $ (\alpha_i\rho_i), (\alpha_i\rho_ie_i), (\alpha_i\rho_iv_i), (\alpha_i\rho_iv_i^2),  p, \dots$
 are linked between themselves by the fact that,  they are *weak limits of respective continuous  functions $\alpha_{i,\epsilon}\rho_{i,\epsilon}, \alpha_{i,\epsilon}\rho_{i,\epsilon}e_{i,\epsilon}, \alpha_{i,\epsilon}\rho_{i,\epsilon}v_{i,\epsilon}, \alpha_{i,\epsilon}\rho_{i,\epsilon}v_{i,\epsilon}^2, p_{i,\epsilon}, \dots$ when $\epsilon\rightarrow 0$. Indeed such products of Radon measures or of Radon measures and $L^\infty$ functions do not make sense mathematically in general. Equations (1-3) are satisfied as a sum of a few distributions obtained by *weak compactness whose sum is the null distribution: the individual products inside the notation of these terms only recall the origin of these terms from a weakly convergent sequence of approximate solutions.\\

This definition of these Radon measures  by passage to the limit of corresponding smooth objects when $\epsilon\rightarrow 0$ appears justified from physics since the molecular structure of fluids forbids that the value $dx=\epsilon$ used to state the conservation laws could be too small:  the cells $\Pi_{1\leq j\leq 3}[x_j-\frac{\epsilon}{2},x_j+\frac{\epsilon}{2}]$ should contain a certain amount of molecules so that the statement of the conservation laws could make sense physically. Therefore the ODEs (6-11) model physics for $\epsilon>0$ very small, but not arbitrarily small, and the limit case $\epsilon=0$ obtained from compactness is only an approximation of the physical situation. This justifies that the above Radon measures are only defined by limits of the corresponding smooth functions and that the apparent products inside them make no sense if not to recall their physical interpretation.\\

In a mathematical viewpoint our method is an inversion of order in the operations of passage to the limit $\epsilon\rightarrow 0$ (in which $\epsilon$ is the side of the small cube $\Pi_{1\leq i\leq n}[x_i-\frac{\epsilon}{2},x_i+\frac{\epsilon}{2}]$ used by physicists to obtain the equations) and resolution of the equations. Indeed to state conservation laws physicists first state (6-8) (with $\frac{d}{dt}w(x,t,\epsilon)$ replaced by $\frac{w(x,t+dt,\epsilon)-w(x,t,\epsilon)}{dt})$, then they pass formally to the limit $\epsilon=0$ to obtain (1-3). Then one attempts to solve directly (1-4). In our method we solve (6-11) for fixed $\epsilon$ then pass to the limit on the solutions to express that $\epsilon>0$ is very small. Finally this method provides Radon measures that satisfy a weak formulation of (1-3).\\

\textbf{6. Taking into account molecular agitation.} To simplify the notation let us consider the case of a single fluid and formula (1) in the familiar form 
\begin{equation} \partial_t\rho +\partial_x(\rho v)=0.\end{equation}
Then formulas (5,6) give
\begin{equation}\frac{d}{dt}\rho(x,t,\epsilon)=\frac{1}{\epsilon}[\rho v^+(x-\epsilon,t,\epsilon)-\rho v^+(x,t,\epsilon)-\rho v^-(x,t,\epsilon)+\rho v^-(x+\epsilon,t,\epsilon)]+\kappa_1(\epsilon).\end{equation}

One considers the 3 cells $\mathcal{C}_{y}=]y-\frac{\epsilon}{2},y+\frac{\epsilon}{2}[ $ with  $y=x,x-\epsilon$ and $x+\epsilon$ with respective constant values $\rho(y,t,\epsilon)$ and $v(y,t,\epsilon)$ inside each cell. Formula (48) describes the transport of the matter around the cell $\mathcal{C}_{x}$: if one considers the cell interface $x-\frac{\epsilon}{2}$ then in the time interval $[t,t+dt]$ the amount of matter $\rho v^+(x-\epsilon,t,\epsilon) dt$ crosses this interface to the right and  the amount of matter $\rho v^-(x,t,\epsilon) dt$ crosses this interface to the left. Therefore (48) models the transport of matter according to its macroscopic velocity $v$ as described in \cite{ColombeauSiam} p. 1906 provided the  condition $|v|dt<\epsilon$ to eliminate the influence of more remote cells during time $dt$.\\

In various circumstances the mathematical proofs imposed or suggested the replacement of $v^\pm$ by $v^\pm+\mu$ with $\mu>0$ large enough. This was done in \cite{ColombeauZeit} pp. 2586-2588 and in \cite{Abreu} pp. 1208-1211. The physical meaning of this replacement is that between times $t$ and $t+dt$ an additional amount of matter $\mu\rho(x-\epsilon,t,\epsilon) dt$ crosses the interface $x-\frac{\epsilon}{2}$ to the right and an additional amount of matter $\mu\rho(x,t,\epsilon) dt$ crosses the interface $x-\frac{\epsilon}{2}$ to the left. This clearly can be interpreted as a model of molecular agitation with mean velocity value (in all senses: to the right and same to the left) of molecules equal to $\mu$, on which the macroscopic velocity $v$ is superposed.\\

This change of $v^\pm$ into $v^\pm+\mu$ does not modify the proof in this paper. One can even change $v^\pm$ into $v^\pm+\frac{\mu}{\epsilon^\beta}, \ 0<\beta<1$: then one simply notices that in (45) $\|v_i^\pm\|_\infty $ becomes $\|v_i^\pm\|_\infty+\frac{\mu}{\epsilon^\beta} $ which gives a bound $O(\epsilon^{1-2\lambda})+O(\epsilon^{1-\beta})$ for the right hand side of (45). This change could play a role  in a search of a physically admissible solution   since it permits to prove that the method in this paper gives the Kruzhkov entropy solution  \cite{Abreu}. This change also  plays a role numerically by adding some vanishing viscosity, see \cite{ColombeauZeit} pp. 2586,2587. Indeed
it has been noticed that the presence of $\mu>0$ large enough is needed in the numerical tests below: in figures 1 and 2 one uses  $\mu=300$. \\

\textbf{7. Numerical calculation of the solutions.} In this section we calculate numerically the solution in the case of the Toumi shock tube problem \cite{Toumi,Munk,FFMunk,FlattenRoe}. In figure 1 we compare the approximate solutions  obtained from our theoretical construction with the result of a  standard "`transport-averaging-pressure correction"' splitted scheme which is a direct extension of the scheme in \cite{ColombeauJDE} section 7 and is described in  \cite{Colombeaunumerique}. We obtain an exact superposition. In figure 2 we compare the approximate solutions from our theoretical construction with the numerical solution obtained  from standard numerical methods of scientific computing in \cite{Munk} pp. 2615,2617, \cite{FFMunk} p. 437, \cite{FlattenRoe} p. 497, in which these authors  insert into  system (1-4) an additional term to render it hyperbolic to improve its mathematical and numerical properties, \cite{Munk} p. 2595. We observe close results in pressure, gas temperature and liquid temperature (same step values and same location of the discontinuities). In gas volume fraction, gas velocity and liquid velocity we still observe same step values and same location of discontinuities, but in these three cases we observe a significant difference in form of peaks: a "`down-peak"' in gas volume fraction and two very neat "`up-peaks"' of well defined height in gas and liquid velocities. The features of these peaks are well defined: they are invariant under changes of the time steps and the CFL number. They are identical for the two different schemes tested in figure 1. Therefore these peaks are not artefacts of calculation since they are obtained from two completely different numerical numerical methods. Further quite similar peaks in shape and values are observed in engineering codes and observations on the gas kick \cite{Avelar,Avelarthese,Bezeira,Choi,Loiola}.\\

 By testing both absence  and presence  of the additional term introduced by many authors to render the system hyperbolic by the same  standard "`transport-averaging-pressure correction"' order 1 numerical scheme   adapted from  \cite{ColombeauSiam} and \cite {ColombeauJDE} section 7,  see \cite{Colombeaunumerique}, one shows that the additional term is responsible of the disappearance of these peaks: without the additional term the scheme produces the result in figure 1, with the additional term the scheme produces the result in \cite{FFMunk,FlattenRoe,Munk}. The additional term is motivated  to ensure hyperbolicity of the system. However in some cases of  fully nonlinear systems of physics it has been observed in \cite{Keyfitz0}-\cite{KeyfitzSan} that hyperbolicity is not always indispensable to produce a well-posed solution.\\

 From the experimental result reproduced in \cite{Avelar}, figure 4, the peak obtained from the asymptotic solution in this paper for system (1-4) and, also  in absence of the additional term,  by the "`transport-averaging-pressure correction"' scheme of \cite{Colombeaunumerique} appears to be related to the gas kick phenomenon. In figure 3 we observe the evolution of the peak in liquid flow rate according to time. Its velocity, top value and width appear of an order of magnitude completely compatible with observations on the gas kick, see \cite{Avelar,Avelarthese,Bezeira,Choi,Loiola} and articles quoted there. These authors report experimental results and results from numerical codes describing the gas and liquid flow rates ($m^3/s$ or $kg/s$) at the top of the well as a function of time. The results reported in this paper concern solutions of the Riemann problem i.e. they represent the gas and liquid physical variables (in particular velocities) inside the tube at a given time. From them one can at once obtain the flow rates at the end of the tube as a function of time: constant velocity before the peaks arrive at the end of the tube, then  peaks for the flow rates as a function of time when the"' kick"' goes out of the tube, completely similar to the results reported by engineers \cite{Avelar,Avelarthese,Bezeira,Choi,Loiola} and articles quoted there. But the experimental results concern real situations with two coaxial tubes and the gas kick comes from the annular domain between the tubes etc, therefore a detailed comparison is not possible within the scope of this mathematical paper. \\

INSERT FIGURE 1\\

Figure 1. \textit{Comparison of the asymptotic solution constructed in this paper (black +) and the result from the transport-averaging-correction scheme of \cite{Colombeaunumerique} with $\delta=0$ (red, continuous line). One observes a perfect coincidence.}\\

INSERT FIGURE 2\\

Figure 2. \textit{Comparison of the asymptotic solution constructed in this paper (black,+) with the numerical solution constructed by various  authors (red, continuous line) by adding an additional term in the equations. This last solution is  obtained here with the schemes of \cite{Colombeaunumerique} after having checked it gives same results.}\\

INSERT FIGURE 3\\

Figure 3. \textit{The liquid flow rate ($kg/s$) at time t=0.03 (top left panel), t=0.06 (top right panel), t=0.12 (bottom panel, in a twice longer tube). One observes that the top value increases very slowly with time, and that its width increases proportionally to time: one interval=2.5 centimeter in each panel.}\\

In the three figures the system of ODEs (6-8) is solved by the explicit Euler order 1 scheme. Space which is 100 meters long is divided into  4000 cells (8000 cells in the bottom panel of figure 3 where we consider a 200 meters long tube); the solution is given at time $t=0. 06,r= \frac{\Delta t}{\Delta x}=0.0005 $ in all tests, $\mu=300$ in figures 1 and 2 and $\mu=$200, 300 and 400 in figure 3 from top left panel to bottom panel with various values of time (t=0.03,0.06,0.12)$, \lambda=1$ in (11) with $\epsilon$ equal to the space step $\Delta x$. \\


\textit{Remark: the case of one fluid.} One has checked that the  adaptation of the  method in  this paper to the case of one fluid has exactly given numerically the known solutions on the four 1-D Toro tests in \cite{Toro1} and (easy extension of the proof to $n$-D, $n=2,3,\dots$ as explained in \cite{ColombeauZeit}) the six 2-D Lax tests in \cite{Lax1,Lax2}.\\

	\textbf{8. Conclusion.} Since various  numerical schemes giving same result for conservative systems can give really different results for nonconservative systems such as (1-4) the search of mathematically well defined solutions is particularly important. In this paper we have proved that one can obtain approximate and Radon measure solutions in a weak sense. By reducing  the system of partial differential equations in the case of these approximate solutions to a system of ordinary differential equations we observe numerically on the standard Toumi shock tube problem that the Radon measures from our method agree  with the numerical solutions previously obtained by other authors  with various different numerical methods, modulo very neat peaks of well defined limited height and width in liquid and gas velocities, which suggest the gas kick phenomenon that appear in our method, engineering codes and experimental observations. Indeed in a subsequent numerical paper, we observe exactly the numerical solution given by the theoretical mathematical proof presented in the present paper.  \\


\centering 
\includegraphics{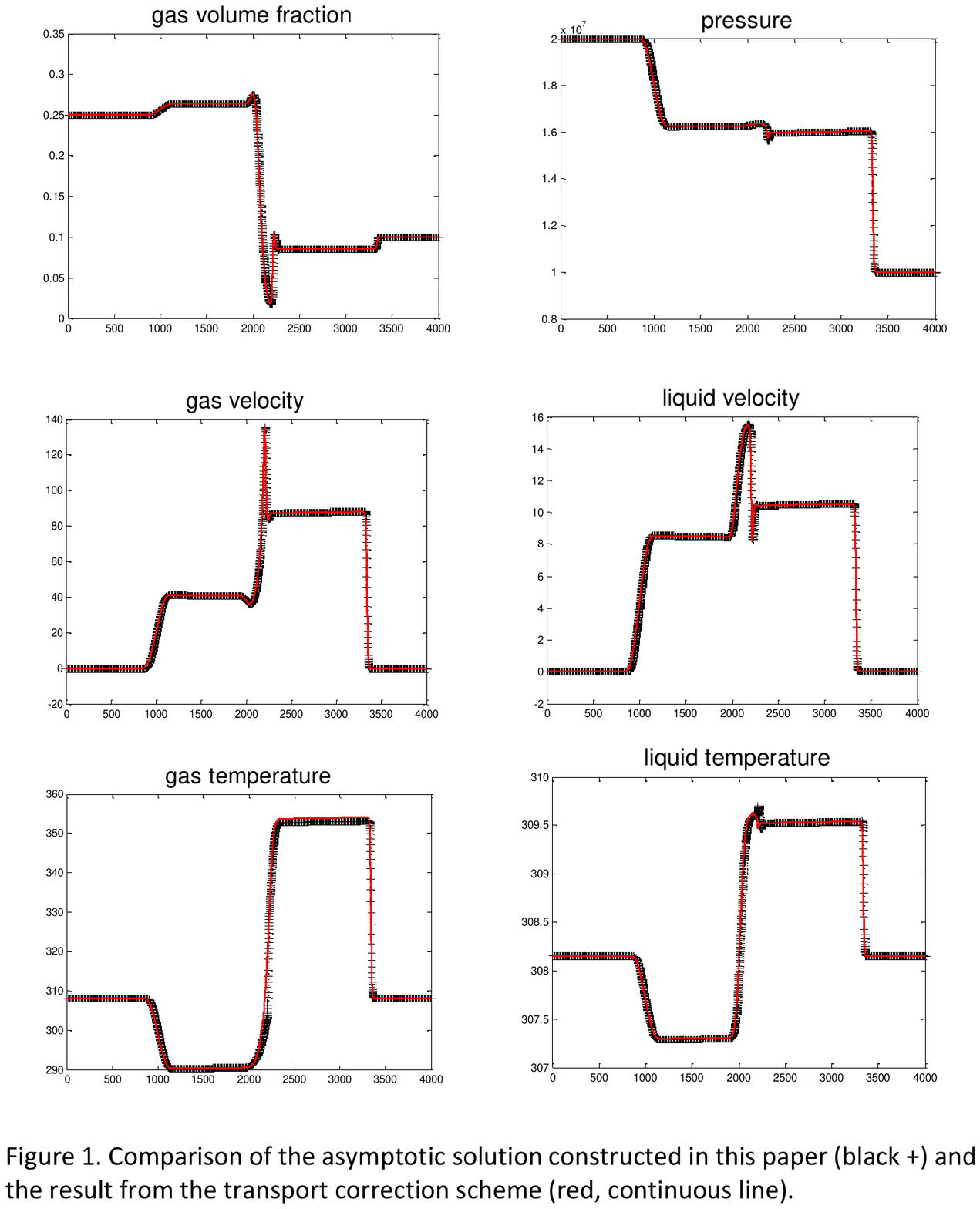}

\centering 
\includegraphics{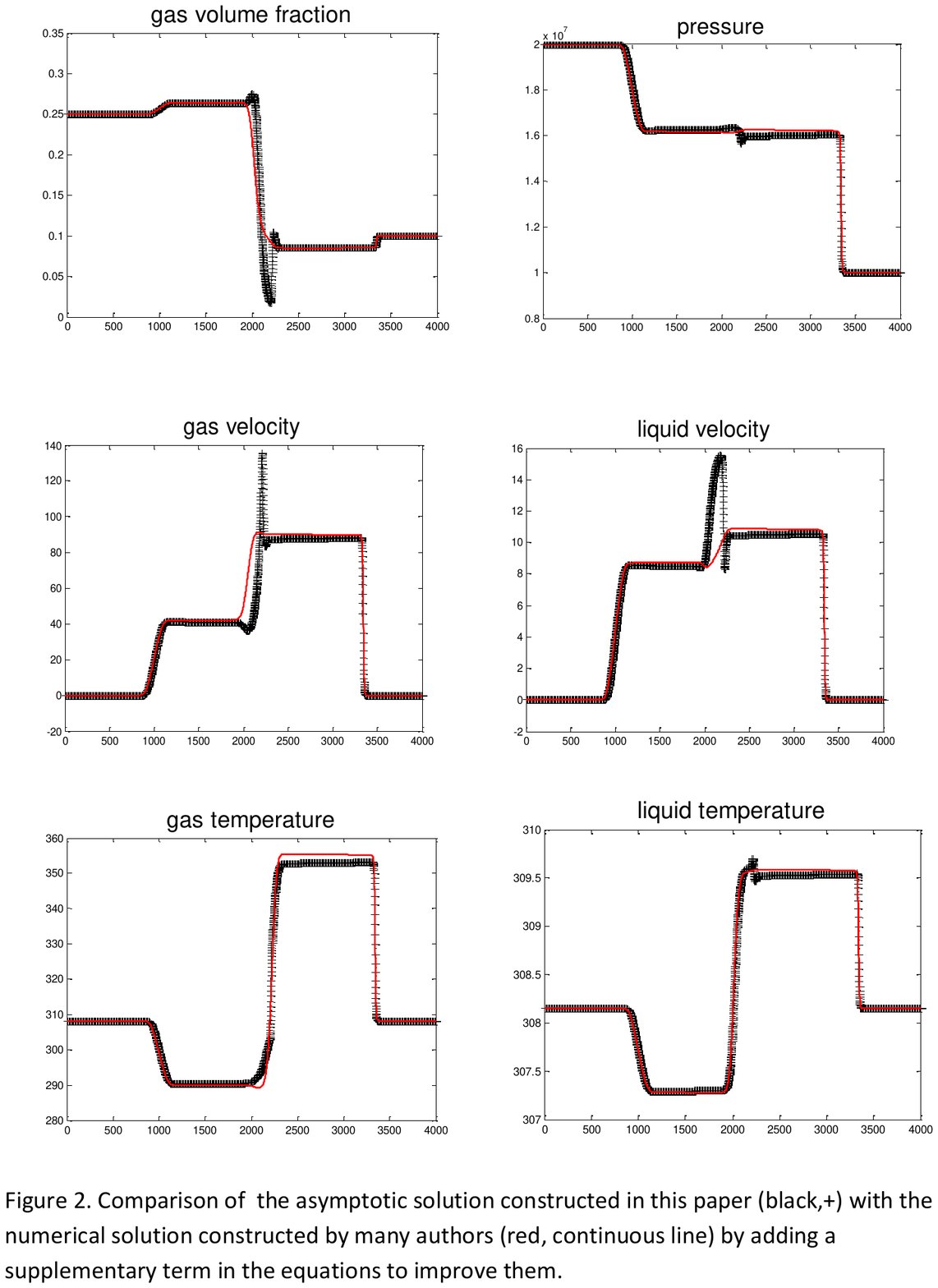}

\centering 
\includegraphics{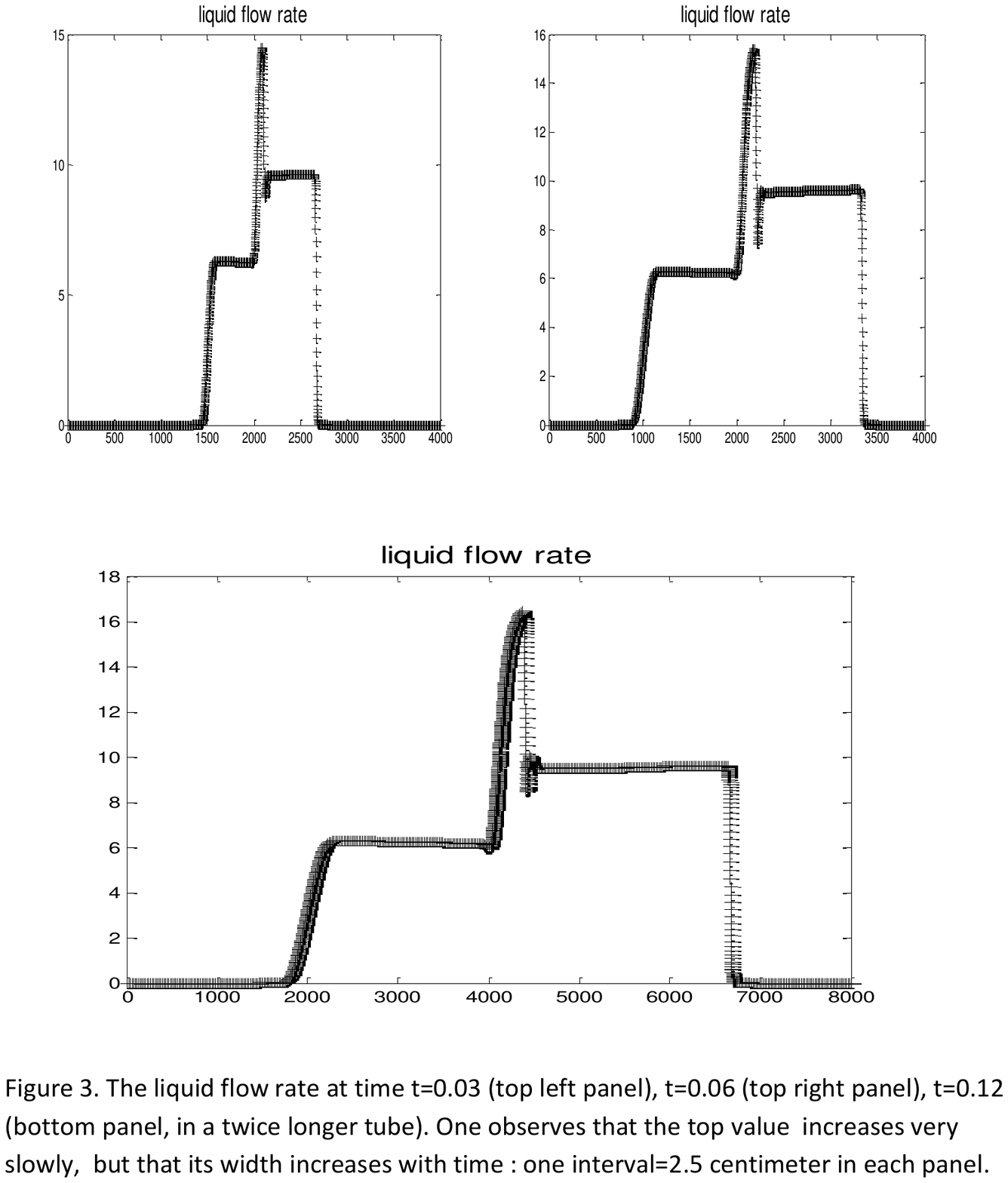}


\begin{thebibliography}{<50>}

\bibitem{Abreu} E. Abreu, M. Colombeau, E. Panov. Weak asymptotic methods for scalar equations and systems. J. Math. Anal. Appl. 44,2, 2016, pp. 1203-1232.



\bibitem{Avelar} C. S. Avelar, P. R. Ribeiro, K. Seperhrnoori. Deep water gas kick simulation. J. Pet. Sci. Eng. 67, 2009, pp. 13-22.

\bibitem{Avelarthese} C. S. Avelar. Well contol modelling: a finite difference approach (portuguese). Master Thesis, State University of Campinas, 26/08/2008.

\bibitem{Bendiksen} K.H. Bendiksen, D. Maines, R. Moe, S. Nuland. The dynamic two-fluid model OLGA. Theory and Application. SPE Prod. Eng. 6, 1991, pp. 171-180.

\bibitem{Bestion} D. Bestion. The physical closure laws in the CATHARE code. Nucl. Eng. Des. 124, 1990, pp. 229-245.

\bibitem{Bezeira} E. M. Bezeira. Study of well control consideruing phase behabior of gas-liquid mixture (portuguese). Master Thesis, State University of Campinas, 07/04/2006.
\bibitem{Brezis} H. Brezis. Functional Analysis, Sobolev spaces and Partial Differential Equations. Springer, New York, Dordrecht, Heidelberg, London. 2010.

\bibitem{Choi} K. I. Choi. A dynamic lagrange-euler numerical model for biphasic flows in oil wells. Master Thesis, CEPETRO, 01/03/1996.



\bibitem{ColombeauSiam} M. Colombeau. A method of projection of delta waves in a Godunov scheme and application to pressureless fluid dynamics. SIAM J. Numer. Anal. 48, 5, 2010, pp. 1900-1919.




\bibitem{ColombeauApplicable} M. Colombeau. Irregular shock wave solutions as continuations of the analytic solutions. Appl. Anal. 94,9,2015, pp.1800-1820.
 
 \bibitem{ColombeauJMP} M. Colombeau. Weak asymptotic methods for 3-D selfgravitating pressureless fluids. Application to the creation and evolution of solar systems from the fully nonlinear Euler-Poisson equations. J. of Mathematical Physics 56, 061506, 20 pages, 2015.

\bibitem{ColombeauZeit} M. Colombeau. Approximate solutions to the initial value problem for some compressible flows. Zeitschrift fur Angewandte Mathematik und Physik. 66, 2015,5, pp. 2575-2599.


\bibitem{ColombeauJDE} M. Colombeau. Asymptotic study of the initial value problem to a standard one pressure model of multifluid flows in nondivergence form. J. Diff. Eq. 260, 2016,1, pp. 197-217.

\bibitem{Colombeaunumerique} M. Colombeau. A numerical approximation of the standard one pressure system of two fluid flows with energy equations under its four versions. arXiv.org. 1808.08467.














\bibitem{EvjeFlatten} S. Evje, T. Flatten. Hybrid Flux-splitting Schemes for a common two fluid model. J. Comput. Physics 192, 2003, pp. 175-210.

\bibitem{Joseph1} K.T. Joseph. Boundary layers in approximate solutions. Trans. Amer. Math. Soc. 314, 1989, pp. 709-726.
\bibitem{Joseph3} K.T. Joseph. P.L. Sadchev. Exact solutions for some nonconservative hyperbolic systems. Intern. J. Nonlinear Mechanics 38,2003, 9, pp. 1377-1386.
\bibitem{Joseph4} K.T. Joseph. Generalized solutions to a Cauchy problem for a nonconservative hyperbolic system. J. Math. Anal. Appl. 2007, 1997, pp.361-387.
\bibitem{Joseph5} K.T. Joseph, Ph LeFloch. Singular limits for the Riemann problem CR Math Acad. Sci. Paris 344, 2007,1, pp.59-64. Analytical Approaches to mutidimensional balance laws. pp. 143-171. Nova Sci. Publi. New York 2006.
\bibitem{Joseph6} K.T. Joseph, M.R. Sahoo. Vanishing viscosity appoach to a system of conservation laws admitting $\delta"$-waves. Comm. Pure Appl. Analysis 12, 2013, 5, pp.2091-2118. 
\bibitem{Joseph7} K.T. Joseph. Exact solution of a nonconservative system in elastodynamics. Electron. J. Diff. Eq. 2015, art. 259, 7p.
 





\bibitem{Keyfitz0} B.L. Keyfitz. Properties of Nonhyperbolic Models for Incompressible Two Phase Flows. Report 77204-3476 of University of Houston, Fields Institute and University of Ohio.


\bibitem{Keyfitzlivro} B.L. Keyfitz. Hold that light! Modelling of traffic flow by Differential Equations. In "`Six Themes on variation"', Robert Hardt editor, Student Mathematical Library, vol. 26, AMS, Providence, Rhode Island, USA, 2004, pp. 127-153.

\bibitem{Keyfitzart1}  B.L. Keyfitz. Admissibility conditions for shocks in conservation laws that change type. Siam J. Math. Anal. 22, 5, 1991, pp. 1284-1292.

\bibitem{Keyfitzart2}  B.L. Keyfitz. Change of type in simple models for two phase flows . In Viscous Profiles and Numerical Methods for shock waves"' ed. M. Shearer, 1991, pp 84-104, SIAM, Philadelphy.

\bibitem{Keyfitzart3}  B.L. Keyfitz. Change of type in three phase flows: a simple analog. J. Diff. Eq. 80, 1989, pp. 280-305.



\bibitem{Keyfitzrev1}  B.L. Keyfitz. Conservation laws that change type and porous medium flow:a review. In "`Modeling and Analysis of Diffusive and Advective Processes in Geosciences"'pp. 122-145.   eds W.E. Fitzgibbon and M.F. Wheeler, SIAM, Philadelphia, 1992.


\bibitem{Keyfitzrev2}  B.L. Keyfitz. Multiphase saturation equations, change of type and inaccessible regions. In Proceedings of the Oberwolfach conference on porous media. pp. 103-116. eds J. Douglas, C.J. Van Duijn and U. Hornung. Birhkhauser, 1993.



\bibitem{Keyfitzrev3}  B.L. Keyfitz. A geometric theory of conservation laws which change type. Zeitschrift fur Angewandte Mathematik und Mechanik, 75, pp. 571-581, 1995.






\bibitem{KeyfitzLopes}  B.L. Keyfitz, M. Lopes-Filho. A geometric study of shocks in equations that change type. J. of Dyn and Diff. Eqs, 6,3, 1994, pp. 351-393.

\bibitem{KeyfitzSan}  B.L. Keyfitz, R. Sanders, M. Sever. Lack of hyperbolicity in the two-fluid model for two-phase incompressible flows. Discrete and continuous dynamical systems, serie B,3,4,2003,pp. 541-563.


\bibitem{Kunzinger2} M. Kunzinger, G. Rein, R. Steinbauer, G. Teschl. Global weak solution of the relativistic Vlassov-Klein Gordon system. Comm. Math. Pkys. 238, 2003, 1-2, pp. 367-378.




\bibitem{Larsen} M. Larsen, E. Hustvedt, P. Hedne, T. Straume. Petra: a novel computer code for simulation of slug flow.  In SPI Annual Technical Conference and Exhibition. SPI 38841, 1997, pp. 1-12.


\bibitem{Lax1} P.D. Lax. Computational Fluid Dynamics. J. Sci. Comput. 31, 2007, pp. 185-193.

\bibitem{Lax2} P.D. Lax. Mathematics and Physics. Bull. A.M.S. 45, 1, 2008, pp. 135-152.

\bibitem{Loiola} C. H. O. Loiola. Compositional well control simulation (portuguese). Master Thesis, CEPETRO, State university of Campinas,22/09/2015. 


\bibitem{FFMunk} P.J. Martinez-Ferrer, T. Flatten, S.T. Munkejord. On the effect of temperature and velocity relaxation in two-phase flow models M2AN 46, 2012, pp. 411-442.

\bibitem{FlattenRoe} A. Morin,T. Flatten, S.V. Munkejord. A Roe scheme for a compressible six-equation two fluid model. International Journal for Numerical Methods in fluids, 2013,72, pp. 478-504.  

\bibitem{Munk} S. T. Munkejord, S. Evje, T. Flatten. A Musta scheme for a nonconservative two-fluid model. SIAM J. Sci. Comput. 31, 4, 2009, pp. 2587-2622.









\bibitem{Toro1} E.Toro. Riemann Solvers and Numerical Methods for Fluid Dynamics. Springer Verlag, 1999.


\bibitem{Toumi} I. Toumi, A. Kumbaro. An approximate linearized  Riemann solver for a two fluid model. J. Comput. Physics 66, 1986, pp. 62-82.

\bibitem{WAHA3} WAHA3 Code Manual. JSI Report IJS-DP-8841, Josef Stefan Institute, Ljubljana, Slovenia, 2004.

\bibitem {Wendroff} H. B. Stewart, B. Wendroff. Two-phase flows: Models and Methods. J. Comput. Phys. 56, 1984, pp. 363-409.




 \end{thebibliography}
 \end{document}